\documentclass[12pt]{article}
\usepackage{amsfonts,amsmath}
\usepackage{amssymb,latexsym}
\usepackage[dvips]{epsfig}
\usepackage{amscd}
 
\title{Triply-graded link homology and Hochschild homology of Soergel bimodules}
\author{ Mikhail Khovanov}
\date{October 14, 2005}
 
\newtheorem{prop}{Proposition}
\newtheorem{theorem}{Theorem}
\newtheorem{lemma}{Lemma}
\newcommand{\oplusop}[1]{{\mathop{\oplus}\limits_{#1}}}
 
\begin{document} 

\maketitle
\baselineskip 14pt
 
\def\R{\mathbb R}
\def\Q{\mathbb Q}
\def\Z{\mathbb Z}
\def\N{\mathbb N} 
\def\C{\mathbb C}
\def\l{\lbrace}
\def\r{\rbrace}
\def\o{\otimes}
\def\lra{\longrightarrow}
\def\Hom{\mathrm{Hom}}
\def\RHom{\mathrm{RHom}}
\def\Id{\mathrm{Id}}
\newcommand{\define}{\stackrel{\mbox{\scriptsize{def}}}{=}}
\def\drawing#1{\begin{center}\epsfig{file=#1}\end{center}}
 \def\yesnocases#1#2#3#4{\left\{
\begin{array}{ll} #1 & #2 \\ #3 & #4
\end{array} \right. }

\begin{abstract} We trade matrix factorizations and Koszul 
complexes for Hochschild homology of Soergel bimodules to 
modify the construction of triply-graded link homology and 
relate it to Kazhdan-Lusztig theory. 
\end{abstract}

\vspace{0.2in} 

{\bf Hochschild homology.} 

\noindent 
Let $R$ be a $k$-algebra, where $k$ is a field,
$R^e=R\otimes_k R^{op}$ be the enveloping algebra of $R,$ and 
$M$ be an $R$-bimodule (equivalently, a left $R^e$-module). 
The functor of $R$-coinvariants 
associates to $M$ the factorspace $M_R=M/[R,M],$ where 
$[R,M]$ is the subspace of $M$ spanned by vectors of the form 
$rm-mr.$ We have $M_R= R\otimes_{R^e} M.$ The $R$-coinvariants 
functor is right exact and its $i$-th derived functor takes $M$ to   
$\mathrm{Tor}_i^{R^e}(R,M).$  The latter space is also denoted 
$\mathrm{HH}_i(R,M)$ and called the $i$-th Hochschild homology of $M.$ 
The Hochschild homology of $M$ is the direct sum 
$$\mathrm{HH}(R,M) \define \oplusop{i\ge 0} \mathrm{HH}_i(R,M).$$ 

\vspace{0.1in} 

To compute Hochschild homology, we choose a resolution of 
the $R$-bimodule $R$ by projective $R$-bimodules
$$  \dots \lra P_2 \lra P_1 \lra P_0 \lra R \lra 0, $$ 
and tensor the resolution with $M:$ 
$$  \dots \lra P_2\otimes_{R^e} M \lra P_1\otimes_{R^e} M 
 \lra P_0\otimes_{R^e} M \lra 0.$$ 
Homology of this complex is isomorphic to the Hochschild homology 
of $M.$ 

Hochschild homology is a covariant functor from the category 
of $R$-bimodules to the category of $\Z_+$-graded $k$-vector spaces. 
In particular, a homomorphism of $R$-bimodules induces a map on their 
Hochschild homology. 
If $R$ has a grading, the Hochschild homology $\mathrm{HH}(M,R)$ 
of a graded $R$-bimodule $M$ is bigraded. 

\vspace{0.1in} 

Any $k$-algebra $R$ has the standard ``bar'' resolution by free 
$R$-bimodules. The polynomial algebra $R=k[y_1, \dots, y_m]$ admits 
a much smaller ``Koszul'' resolution by free $R^e=R\otimes R$-modules 
(since $R$ is commutative,  $R= R^{op}$), given by the tensor product 
(over $R^e$) of complexes 
$$0 \lra R^e \xrightarrow{y_i \otimes 1 - 1\otimes y_i} R^e \lra 0 $$ 
for $i=1, \dots, m.$ 
The resolution has length $m$ (spanning homological degrees between $0$ and 
$m$) and its total space is naturally the tensor 
product of $R^e$ and the exterior algebra on $m$ generators. 

Thus, the Hochschild homology of a bimodule $M$ over the polynomial 
algebra $R$ is the homology of the complex built out of $2^m$ copies 
of $M$ (the $i$-th term of the complex consists of  $m$ choose $i$   
copies of $M$), with the differential built out of multiplications by 
$y_i\otimes 1 - 1\otimes y_i.$ More precisely, denote by $C(M)$ the 
chain complex 
$$ 0 \lra C_m(M) \lra \dots \lra C_1(M) \lra C_0(M) \lra 0$$ 
where
 $$C_j(M) = \oplusop{I\subset \{1, \dots, m\}, |I|=j} M\otimes_{\Z} \Z 
 \lfloor I \rfloor ,$$
the sum over all subsets $I$  of cardinality $j$ and $\Z \lfloor I\rfloor$ 
being rank 1 free abelian group generated by the symbol $\lfloor I\rfloor.$ 
The differential has the form 
$$ d(m\otimes \lfloor I\rfloor) = \sum_{i\in I} \pm (y_i m - my_i) 
\otimes \lfloor I\setminus \{i\}  \rfloor, $$
where we choose the minus sign if $I$ contains an odd number of elements less 
than $i.$ Then $\mathrm{HH}(R,M) \cong \mathrm{H}(C(M)).$
Clearly, 
$$\mathrm{HH}_m(R,M) = \{m\in M| rm=mr \hspace{0.1in} 
\forall r\in R\} =M^R,$$ 
the $R$-invariants subspace of $M.$ 

The Hochschild \emph{cohomology} groups 
of $M$ are the derived functors of the left exact 
functor of $R$-invariants, evaluated on $M$.  For the polynomial 
algebra $R,$ Hochschild homology and cohomology are isomorphic,  
$$ \mathrm{HH}_i(R,M) \cong \mathrm{HH}^{m-i}(R,M)$$ 
for any bimodule $M.$ This property, which can be explained by the self-duality of 
the Koszul resolution of $R,$ does not extend to arbitrary algebras. It implies that 
in all constructions described below we can substitute 
Hochschild cohomology for Hochschild homology without any gain or loss. 

Hochschild cohomology of $R$-bimodules $M,$  for any $R,$ are 
covariant (rather than contravariant) in $M,$ just like Hochschild homology. 
For a thorough treatment of Hochschild (co)homology we refer the reader to 
Loday's book [L], and to Kassel [Ka] for a very brief introduction. 

\begin{figure} [htb] \drawing{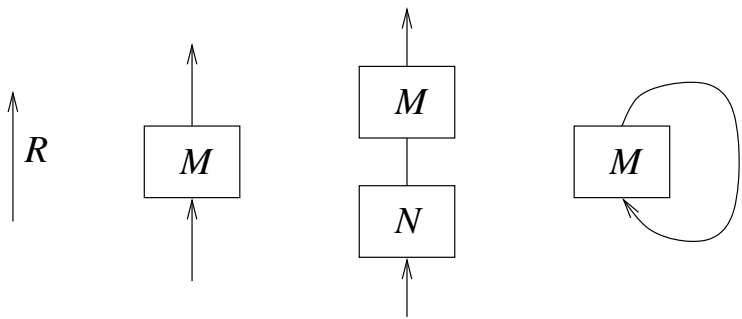} \caption{From left to right: 
$R$-bimodules $R, M, N\otimes_R M;$ Hochschild homology of $M.$}\label{bimodules} 
\end{figure} 

\vspace{0.1in} 

For our purposes the figure~\ref{bimodules} diagrammatical calculus will come in 
handy. Depict 
the $R$-bimodule $R$ by an oriented line, an $R$-bimodule $M$ by a box extended 
by two lines symbolizing the left and right actions of $R.$ Concatenation of boxes 
is interpreted as the tensor product $N\otimes_R M,$ and the closure of two ends 
of boxed $M$ as taking the Hochschild homology of $M.$

\vspace{0.2in} 

{\bf Soergel bimodules.} 

\noindent 
Let $R'=\Q[x_1, \dots, x_m]$ and $R_i'=\Q[x_1, \dots, x_{i-1}, x_i+x_{i+1}, 
x_i x_{i+1}, x_{i+2}, \dots, x_m].$ The ring $R_i'$ is a subring of $R'$ of 
polynomials which are symmetric in $x_i$ and $x_{i+1}.$ Equivalently, 
$R_i'$ consists of polynomials invariant under the action of the symmetric 
group $S_2$ which exchanges $x_i$ and $x_{i+1}.$ The ring $R'$ is 
a free $R_i'$-module of rank $2.$ 

Introduce a grading on $R'$ and $R_i'$ by placing each $x_i$ 
in degree two. Then $R' \cong R_i'\{2\} \oplus R_i'$ as a graded 
$R_i'$-module, where $\{2\}$ is the grading shift up by two.  

Let $B_i'= R' \otimes_{R'_i} R'.$ It's a graded $R'$-bimodule, 
free of rank two as a left $R'$-module and as a right $R'$-module. 

\vspace{0.1in} 

We will call $B_i',$ their tensor products, and other related bimodules 
\emph{Soergel bimodules}, after Wolfgang Soergel, who introduced them 
and explained their importance for the infinite-dimensional 
representation theory of simple Lie algebras and (closely related) 
Kazhdan-Lusztig theory [S1], [S2]. 

Consider degree-preserving homomorphisms of graded $R'$-bimodules:  
$$ \mathrm{br}'_i: B_i'\lra R', \hspace{0.2in} 
     \mathrm{rb}'_i: R'\{2\} \lra B_i', $$ 
where $\mathrm{br}'_i(1\otimes 1) = 1$ and 
$$ \mathrm{rb}'_i(1) = (x_i-x_{i+1})\otimes 1 + 1\otimes 
 (x_i-x_{i+1}).$$ 

We would like to consider arbitrary tensor products (over $R'$)  
$$B_{i_1}'\otimes B_{i_2}'\otimes \dots \otimes B_{i_n}'$$ 
Bimodules $B_i'$ and all their tensor products have a trivial 
direction, in the following sense. Let 
$$R= \Q[x_1-x_2, \dots, x_{m-1}-x_m]$$  
be the polynomial ring generated by consecutive differences 
of variables $x_1, \dots, x_n$ (note that in [KR2]  $R$ denotes 
another ring).  This is a subring of $R'$ and 
we can write $ R'= R\otimes_{\Q} \Q[x_j],$ for any $j$ (we 
choose $j=1$ from now on). The 
permutation action $x_i\leftrightarrow x_{i+1}$  on $R'$ restricts 
to that on $R$ and we define $R_i\subset R$ to be the ring 
of $S_2$-invariants. Further, set 
$$B_i = R\otimes_{R_i} R $$ 
We have 
$$R' \cong R\otimes \Q[x_1], \hspace{0.2in} R_i'\cong R_i\otimes \Q[x_1], 
\hspace{0.2in} B_i' \cong  B_i \otimes \Q[x_1],$$ 
and bimodule homomorphisms $\mathrm{br}'_i , \mathrm{rb}'_i$ 
restrict to bimodule homomorphisms, denoted $\mathrm{br}_i $ and $\mathrm{rb}_i,$ 
between $R$-bimodules $B_i$ and $R.$ 

\emph{Example:} $m=2.$ Let $y=x_1-x_2.$ Then $R=\Q[y], 
$ $R_1= \Q[y^2]$ and $B_1= \Q[y]\otimes_{\Q[y^2]} \Q[y].$ 
The bimodule homomorphism $\mathrm{rb}_1$ is defined by 
$\mathrm{rb}_1(1) = y \otimes 1 + 1\otimes y.$   

We have an $R'$-bimodule isomorphism, for an arbitrary tensor product: 
$$ B'_{j_1}\otimes_{R'} B'_{j_2}\otimes \dots \otimes_{R'} B'_{j_n} \cong 
  (B_{j_1}\otimes_R B_{j_2} \otimes_R \dots \otimes_R B_{j_n} ) \otimes_{\Q} 
   \Q[x_1].$$ 
Thus, the left hand side can be recovered from the tensor product 
$B_{j_1}\otimes_R B_{j_2} \otimes_R \dots \otimes_R B_{j_n}$ and 
vice versa.

\vspace{0.35in} 

\noindent
{\bf Soergel bimodules and a braid group action (after Rapha\"el Rouquier).} 

\noindent 
Rapha\"el Rouquier [R] pointed out and explored an explicit relation between 
Soergel 
bimodules $B_i$ and the braid group. We recall his results, taking the liberty to use
our conventions. Assign to the braid generator $\sigma_i$ the complex 
$F(\sigma_i)$ of graded $R$-bimodules 
$$F(\sigma_i): \hspace{0.3in} 
  0 \lra R\{2\} \stackrel{\mathrm{rb}_i}{\lra} B_i \lra 0,$$ 
with $B_i$ placed in cohomological degree $0.$ 
To the braid generator $\sigma_i^{-1}$ assign the complex 
$F(\sigma_i^{-1})$ of graded $R$-bimodules 
$$F(\sigma_i^{-1}): \hspace{0.3in}
  0 \lra B_i\{-2\} \stackrel{\mathrm{br}_i}{\lra} R\{-2\} \lra 0,$$
with $B_i\{-2\}$ placed in cohomological degree $0$ (we suggest that 
the reader compares these complexes with figure 6 in [KR2] which 
assigns certain complexes of matrix factorizations to braid generators.) 

To a braid word 
\begin{equation}\label{sigma}
 \sigma = \sigma_{j_1}^{\epsilon_1} \sigma_{j_2}^{\epsilon_2} 
 \dots \sigma_{j_n}^{\epsilon_n}, \hspace{0.2in} \epsilon_i\in \{1, -1\},
\end{equation}  
assign the tensor product (over $R$) of the above complexes and denote 
it by $F(\sigma),$  
For instance, to $\sigma_2 \sigma_3^{-1}\sigma_1$ we assign 
the complex of bimodules 
$$ F(\sigma_2)\otimes_R F(\sigma_3^{-1}) \otimes_R F(\sigma_1).$$ 
We consider the category $\mathcal{B}(R)$ of complexes of graded $R$-bimodules 
up to chain homotopies and view $F(\sigma)$ as an object of $\mathcal{B}(R).$ 

\begin{prop} If braid words $\sigma$ and $\widetilde{\sigma}$ represent 
the same element of the braid group then complexes $F(\sigma)$ and 
$F(\widetilde{\sigma})$ are isomorphic in $\mathcal{B}(R).$ 
\end{prop} 

See Rouquier [R, Section 3] for a proof of this proposition and of more general 
results. In particular, the tensor product $F(\sigma_i)\otimes_R F(\sigma_i^{-1})$ is 
chain homotopy equivalent to the complex $0 \lra R \lra 0$ of $R$-bimodules. 

The tensor product over $R$ is a bifunctor 
$$ \mathcal{B}(R) \times \mathcal{B}(R) \lra \mathcal{B}(R),$$ 
and each object $N$ of $\mathcal{B}(R)$ gives rise to an endofunctor 
of the category $\mathcal{B}(R)$ which takes a complex $M$ 
to the tensor product $M\otimes_R N.$ 

The above proposition says that bimodule complexes $F(\sigma_i)$ 
give rise to a (weak) braid group action on $\mathcal{F}.$ Rouquier 
shows that the action is "genuine", i.e. comes with a transitive 
system of isomorphisms [D]. 

\vspace{0.2in} 

{\bf Graphical presentation.} 

\begin{figure} \drawing{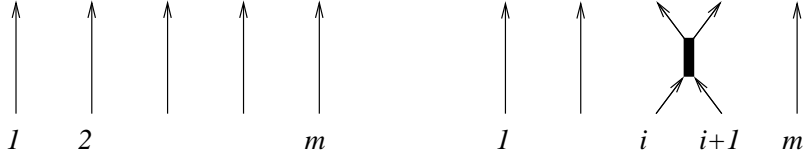} \caption{To left and right pictures we 
assign $R$-bimodules $R$ and $B_i,$ respectively.} \label{lines}
\end{figure}

\noindent 
We now refine the figure~\ref{bimodules} diagrammatics to fit our
situation. 
To $m$ parallel oriented vertical lines we assign the $R'$-bimodule $R'.$ 
The lines symbolize generators $x_1, \dots, x_m$ of $R'.$ The same 
notation will be used to depict $R.$ Bimodules $B_i'$ and $B_i$ will be
assigned to a diagram with a wide edge as in [KR2] bounded by four 
oriented lines with endpoints in the $i$-th an $(i+1)$-st positions, the 
rest of the diagram consisting of oriented lines, see figure~\ref{lines} right. 

\begin{figure} \drawing{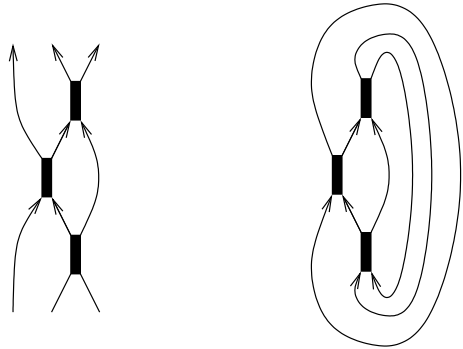}\caption{To the 
left picture we assign  $R_3$-bimodule $B_2\otimes B_1 \otimes B_2$; 
to the right picture--the Hochschild homology of this bimodule.} 
\label{product} 
\end{figure}

To a composition of diagrams with wide edges we assign the tensor 
product of corresponding bimodules, and to the closure of a composition -- 
the Hochschild homology of the tensor product, see figure~\ref{product}.  

\vspace{0.2in} 

{\bf Link homology.} 

\noindent 
Let $\sigma$ be a braid word (see formula (\ref{sigma})). The Rouquier 
complex $F(\sigma)$:  
$$ \dots \stackrel{\partial}{\lra} F^j(\sigma) \stackrel{\partial}{\lra} 
F^{j+1}(\sigma) \stackrel{\partial}{\lra} \dots $$ 
has $n+1$ nontrivial terms, where $n$ is the length of $\sigma.$ 
Each term $F^j(\sigma)$ is a direct sum of graded bimodules 
which are tensor products of $B_i$'s (tensoring with $R$ doesn't do 
anything to a bimodule). One of the summands, for a suitable $j$, is 
$R,$ which we view as the tensor product of zero number of $B_i$'s. 

The Hochschild homology $\mathrm{HH}(R,F^j(\sigma))$ of 
the bimodule $F^j(\sigma)$ is 
a bigraded $\Q$-vector space. Taking the Hochschild homology of each term, we 
obtain a complex of bigraded vector spaces 
$$ \dots \stackrel{\mathrm{HH}(\partial)}{\lra} 
 \mathrm{HH}(R,F^j(\sigma)) \stackrel{\mathrm{HH}(\partial)}{\lra} 
\mathrm{HH}(R,F^{j+1}(\sigma)) \stackrel{\mathrm{HH}(\partial)}{\lra} \dots $$ 
Its cohomology, which we denote $\mathrm{HHH}(\sigma),$ is a 
triply-graded $\Q$-vector space. 

\begin{theorem} \label{thm} Up to an overall shift in the grading, 
$\mathrm{HHH}(\sigma)$ is an invariant of oriented links 
and, up to isomorphism, depends only on the closure of $\sigma.$ This homology 
theory is isomorphic to the reduced homology $\overline{H}(\sigma)$ 
as defined in [KR2, end of Section 1]. 
\end{theorem} 

The theorem implies that the Euler characteristic of 
$\mathrm{HHH}(\sigma)$ is the HOMFLYPT link polynomial 
[HOMFLY], [PT]. By introducing a fractional 
$\frac{1}{2}\Z$-trigrading and a suitable shift, as in Wu [W], 
the grading indeterminancy can be renormalized away. 

\vspace{0.1in} 

\emph{Sketch of proof.} We assume familiarity with [KR2]. Homology 
groups $\overline{H}(\sigma)$ and $\mathrm{HHH}(\sigma)$ 
have similar definitions. In both cases we resolve each crossing 
of the braid in two ways and obtain $2^n$ resolutions of $\sigma.$ 
Each resolution $D$ is a braid diagram of a planar 
graph which is the closure of a concatenation of wide edges, see 
[KR2] and figures~\ref{lines}, \ref{product} above. 

In this paper we assign to $D$ the Hochschild homology $\mathrm{HH}(R,B(D)),$ 
where $B(D)$ denotes the $R$-bimodule which is the tensor 
product of $B_i$'s over all wide edges of $D.$ In [KR2] to $D$ 
we assigned $CH(D),$ the cohomology of the tensor product of 
matrix factorizations over all wide edges and arcs of $D.$ 
In both cases we finish by arranging these groups, over all 
resolutions $D,$ into a complex and taking its cohomology.
Cohomology groups 
$$H(\sigma) = H(\oplusop{D} CH(D), \partial)$$ 
have a "trivial" variable and can be written as  
$$H(\sigma)= \overline{H}(\sigma) \otimes \Q[x],$$ 
with all the complexity carried by $\overline{H}(\sigma).$

Also, the variable $a$ of [KR2] is nearly superfluous and was used, for the 
most part, to keep track of the grading. Setting $a=0$ 
and then following the construction of [KR2] results in 
link homology which is the direct sum of two copies of $H(\sigma),$ 
with a relative shift in the trigrading. 

When $a=0,$ all matrix factorizations in [KR2] turn into 
Koszul complexes. The factorization associated with a wide edge 
[KR2, figure 2] becomes the Koszul complex of the 
sequence $(x_1 +x_2 -x_3 -x_4, x_1x_2-x_3x_4)$ in 
the polynomial ring $\Q[x_1, x_2, x_3 , x_4].$ The Koszul 
complex of this regular sequence has cohomology only in the rightmost 
degree. The cohomology is the quotient 
$$ \Q[x_1,\dots, x_4]/( x_1 +x_2 -x_3 -x_4, x_1x_2-x_3x_4),$$
naturally isomorphic to the bimodule $B_1'$ over the 
polynomial algebra $R'=\Q[x_1, x_2]$ (see the section ``Soergel bimodules''
above). Multiplication by $x_1$ and $x_2$ corresponds to the right action 
of $R'$ on $B_1',$ multiplication by  $x_3$ and $x_4$--to the left action.   
Equalities $x_1+x_2=x_3+x_4$ and $x_1x_2=x_3x_4$ match the 
definition of $B_1'$ as the tensor product $R'\otimes_{R'_1} R'$ over 
the subalgebra $R'_1$ of symmetric polynomials in $x_1,x_2.$ 

\begin{figure} \drawing{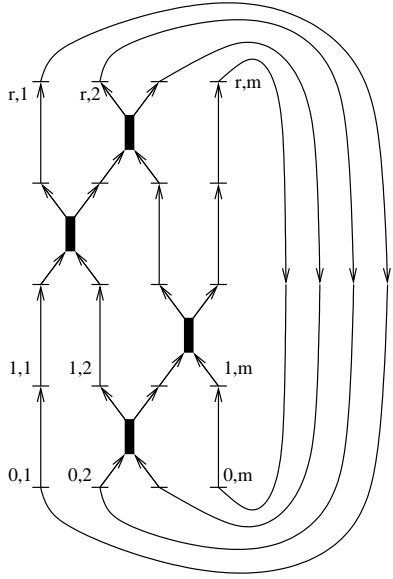} \caption{A resolution $D$ of 
a braid's closure} \label{comp}
\end{figure}

Thus, the relation between [KR2] and Soergel bimodules becomes clear, 
at least locally, for a single wide edge. Globally, a resolution $D$ 
consists of several wide edges and a number of arcs, as shown in 
figure~\ref{comp}. Put marks 
$(i,j)$ on arcs of $D,$ for $0\le i\le r$ and  $1\le j\le m,$ 
where $r$ is the number of wide edges. 

Consider the polynomial 
ring $\widetilde{R}=\Q[x_{i,j}]$ in $(r+1)m$ variables. 
The Koszul complex of $D$ (in the $a=0$ case) corresponds to the 
following sequence of $(r+1)m$ elements of $\widetilde{R}.$ 
For each $1\le i\le r$ we have $m$ elements, two of which 
are $x_{i,s}+x_{i,s+1}-x_{i-1,s}-x_{i-1,s+1},$ 
$ x_{i,s} x_{i,s+1}-x_{i-1,s} x_{i-1,s+1}$ and the rest are 
$x_{i,j}-x_{i-1,j}$ for $1\le j\le m,$ $j\not= s, s+1.$ Here $s$ 
is the position (counting from the left) 
of the $i$-th wide edge of $D$; $s$ is a function of $i.$ 
The remaining $m$ elements of the sequence are 
$x_{0,j}-x_{r,j}$  for $1\le j\le m.$ 

\begin{lemma} The first $rm$ elements on the above list 
constitute a regular sequence in $\widetilde{R}.$ 
\end{lemma} 

To prove the lemma, start with the bottom $m$ variables 
$x_{0,1}, \dots, x_{0,m}$ and work your way up. Each 
time we add a new layer of $m$ variables, we encounter $m$ 
new elements 
$x_{i,s}+x_{i,s+1}-x_{i-1,s}-x_{i-1,s+1},$  
$ x_{i,s} x_{i,s+1}-x_{i-1,s} x_{i-1,s+1},$ and  
$x_{i,j}-x_{i-1,j}$ for $1\le j\le m,$ $j\not= s, s+1.$
These constitute a regular sequence, for they can be matched 
with the  new variables (for any $\Q$-algebra $S$ and 
any $f\in S,$ the element $x^k-f$ is not a zero divisor 
in $S[x]$). $\square$ 

\vspace{0.1in} 

By the time we reach the top layer of variables, we 
take the quotient of $\widetilde{R}$ by the above 
$rm$ variables. The quotient ring is naturally isomorphic 
to the Soergel bimodule $B'(D)$ assigned to $D,$ the 
tensor product of $r$ bimodules $B'_s$ ($s$ is as before, 
and depends on the layer of the diagram). The Koszul 
complex of $D$ is quasiisomorphic to the Koszul 
complex of the quotient ring $B'(D)$ assigned to the sequence of 
the remaining $m$ elements $x_{0,j}-x_{r,j},$ for $1\le j\le m.$  
The latter complex computes the Hochschild homology 
of $B'(D).$ Thus, the Hochschild homology of $B'(D)$ is 
naturally isomorphic to the homology of the Koszul 
complex of $D,$ the latter isomorphic to the direct sum of 
two copies of $CH(D).$ 

Downsizing from $m$ variables $x_{i,1}, \dots, x_{i,m}$ to 
their differences gives us an isomorphism between the 
reduced homology $\overline{CH}(D)$ and 
the Hochchild homology of $B(D).$ These isomorphisms, 
over all $D,$ respect the differentials in the complexes 
computing $\overline{H}(\sigma)$ and $\mathrm{HHH}(\sigma),$ 
leading to the desired isomorphism 
$$\mathrm{HHH}(\sigma)\cong \overline{H}(\sigma).$$  
It was shown in [KR2] that $\overline{H}(\sigma)$ is a 
link invariant, up to a grading shift. Therefore, the same is 
true of $\mathrm{HHH}(\sigma).$ 

$\square$ 

\vspace{0.15in} 

The isomorphism between $\overline{H}$ and $\mathrm{HHH}$ 
is nontrivial on their trigradings. First of all, a not so natural 
shift by $(-1,1,0)$ was built into the definition of $\overline{H}$ due 
to the presence of $a.$ After shifting $\overline{H}$ back by $(1,-1,0),$ 
both homology groups $\overline{H}(\sigma_{\ast})$ and 
$\mathrm{HHH}(\sigma_{\ast})$ of the trivial 1-strand braid 
$\sigma_{\ast}$ become one-dimensional vector spaces sitting in tridegree 
$(0,0,0).$ After this shift is accounted for, the trigradings of the 
two theories relate as follows. 

The third  gradings of $\overline{H}$ and $\mathrm{HHH}$ perfectly match. 
The third grading is the ``cohomological'' grading of both 
theories which comes last into the definition and not visible 
on the homology groups $CH(D)$ and $H(R,B(D))$ 
assigned to resolutions.

The Hochschild grading on $\mathrm{HHH}$ matches the Koszul 
grading on $\overline{H}$ taken with the minus sign, the sign  
due to the difference in conventions.  

Finally, the second grading on $\overline{H}$ equals 
the grading on $\mathrm{HHH}$ by $\deg(x_i)$ minus 
the Hochschild grading, due to the normalization in [KR2, page 2] 
giving $d$ bidegree $(1,1).$ 

\vspace{0.15in} 

{\bf Hecke algebras, Soergel bimodules and Kazhdan-Lusztig theory.} 

\noindent
The Iwahori-Hecke algebra $H_m$ of the symmetric group is a $\Z[q,q^{-1}]$-algebra
 with generators $T_1,\dots T_{m-1}$ and relations 
\begin{eqnarray*} 
  (T_i-q^2)(T_i+1) & = & 0, \\ 
   T_i T_{i+1}T_i & = & T_{i+1} T_i T_{i+1}, \\
   T_i T_j  & = & T_j T_i \hspace{0.2in} \mathrm{for} \hspace{0.2in} 
   |i-j|>1.
 \end{eqnarray*} 
(our $q$ is usually denoted $q^{\frac{1}{2}}$ in the literature).  
$H_m$ is a free $\Z[q,q^{-1}]$-module of rank $n!$ and has 
a natural basis  $\{T_w\}_w,$ over all permutations $w\in S_m,$ where 
$T_w=T_{i_1}\dots T_{i_r}$ if $s_{i_1}\dots s_{i_r}$ is a reduced 
presentation of $w$ as the product of transpositions $s_i=(i,i+1).$ 
We denote $r$ by $l(w)$ and call it the length of $w.$ In our notations, 
$T_{s_i}=T_i.$ Left multiplication 
by $T_i$ in this basis has the form 
$$T_i T_w = \yesnocases{T_{s_iw}}{\mathrm{if} 
  \hspace{0.15in} l(s_iw)> l(w)}{q^2 T_{s_iw} + (q^2-1)T_w}{\mathrm{if} 
  \hspace{0.15in} l(s_iw) < l(w)}$$ 
The involution $\iota$ of $H_m$ takes $q$ to $q^{-1},$ 
 $T_w$ to $(T_{w^{-1}})^{-1}$ and acts on generators by 
$$ \iota(T_i) = q^{-2} T_i +(q^{-2}-1).$$ 
The element $C_i'= q^{-1}(1+T_i)$ is fixed by 
$\iota.$ 

\begin{prop} The Hecke algebra has a basis $\{C_w'\},$ over all $w\in S_m,$ 
with the elements uniquely determined by the conditions 
$\iota(C_w')= C_w'$ and
$$ C_w'= q^{-l(w)} \sum_{y\le w} P_{y,w}(q) T_y,$$ 
where $P_{y,w}(q)\in \Z[q^2], $ $P_{w,w}=1$ and 
$\deg P_{y,w}< l(w) -l(y)$ if $y\not= w.$ The inequality sign 
under the sum sign refers to the Bruhat partial order on $S_m.$ 
\end{prop} 

This proposition is due to Kazhdan and Lusztig [KL1] and admits a direct 
combinatorial proof. The basis $\{C_w'\}$ is called the Kazhdan-Lusztig 
basis (they also introduced a related basis $\{C_w\}$ which won't appear 
in this exposition).  The incredibly hard result, though, is the following. 

\begin{prop} Coefficients of polynomials $P_{y,w}$ are nonnegative 
integers. The multiplication in the basis $\{C_w'\}$ has coefficients 
in $\mathbb{N}[q,q^{-1}].$ 
\end{prop} 

Kazhdan-Lusztig's proof of this result [KL2] is essentially a categorification: 
coefficients of $P_{y,w}$ are realized as dimensions of cohomology 
groups associated to simple perverse sheaves on the flag variety 
$G/B,$ where $G=SL(m,\mathbb{C})$ and $B$ a Borel subgroup of $G.$ 
 Multiplication in the basis $\{C_w'\}$ 
corresponds to the convolution of correspondences
on the product variety $G/B\times G/B.$ Both positivity results 
rely on  Beilinson, Bernstein and Deligne's theory of perverse 
sheaves [BBD] (see also [GM]), including the decomposition theorem, which, in 
turn, requires Deligne's theory of weights and mixed $l$-adic 
sheaves (an outgrowth of Deligne's proof of the Weil conjectures). 
The latter is based on Grothendieck's etale cohomology theory 
of varieties in finite characteristic. 
A characteristic zero alternative approach, via mixed Hodge modules, 
was developed by M.~Saito (see [Sa], [T] and references therein). 
 
\vspace{0.1in} 

Moreover, the work of Beilinson-Bernstein [BB] and 
Brylinski-Kashiwara [BK] on localization relates the whole story 
to infinite-dimensional representations  of $\mathfrak{sl}_m$ 
and implies that $P_{y,w}(1)$ describe multiplicities 
of simple modules in the Verma modules for $\mathfrak{sl}_m.$ 

\vspace{0.1in} 

Here are some additional references on these topics. See [H] for Hecke algebras  
and Kazhdan-Lusztig polynomials, [Sp], [Mi] for localization, 
[M], [Ki] for intersection homology and perverse sheaves, [KW] for 
l-adic and perverse sheaves and relation to KL polynomials. 

\vspace{0.1in} 

To explain how Soergel bimodules fit into the picture, 
we rewrite the defining relations for the Hecke algebra  
via the generators $C_i'=C_{s_i}'.$ The relations become 
\begin{eqnarray*} 
     C_i'^2 & = & (q+q^{-1}) C_i', \\
     C_i' C_{i+1}' C_i' + C_{i+1}' & = & 
     C_{i+1}' C_i' C_{i+1}'  +  C_i', \\
     C_i' C_j' & = & C_j' C_i'   
    \hspace{0.2in} |i-j|>1. 
\end{eqnarray*} 
To match the bimodule $B_i$ with the generator $C_i',$ 
shift its grading down by $1$ and denote by 
$$\mathbf{B}_i = B_i \{ -1\}.$$ 
Now recall one of Soergel's results [S1]. 
\begin{prop} There are isomorphisms of graded $R$-bimodules 
 \begin{eqnarray*} 
    \mathbf{B}_i \otimes_R \mathbf{B}_i & \cong &  \mathbf{B}_i\{1\} 
 \oplus \mathbf{B}_i \{-1\}, \\
    (\mathbf{B}_i \otimes_R \mathbf{B}_{i+1} \otimes_R \mathbf{B}_{i}) 
 \oplus \mathbf{B}_{i+1} 
    & \cong & 
    (\mathbf{B}_{i+1} \otimes_R \mathbf{B}_i \otimes_R \mathbf{B}_{i+1})
 \oplus \mathbf{B}_i, \\   
 \mathbf{B}_i \otimes_R \mathbf{B}_j & \cong & 
 \mathbf{B}_j \otimes_R \mathbf{B}_i \hspace{0.3in} |i-j|>1. 
\end{eqnarray*} 
\end{prop} 
Only the middle isomorphism is non-trivial. Let
$$\mathbf{B}_{i,i+1} = R\otimes_{R_{i,i+1}} R\{-3\},$$
where $R_{i,i+1}$ is the ring of invariants under the action of $S_3$ 
on $R$ permuting $x_i,x_{i+1},x_{i+2}.$ Soergel shows 
\begin{eqnarray*} 
   \mathbf{B}_i \otimes_R \mathbf{B}_{i+1} \otimes_R \mathbf{B}_{i} &  
     \cong & \mathbf{B}_{i,i+1}\oplus \mathbf{B}_i,   \\
  \mathbf{B}_{i+1} \otimes_R \mathbf{B}_i \otimes_R \mathbf{B}_{i+1} &  
     \cong & \mathbf{B}_{i,i+1}\oplus \mathbf{B}_{i+1}, 
\end{eqnarray*} 
which implies the middle isomorphism. 

The above isomorphisms between tensor products of Soergel bimodules 
lift defining relations in the Hecke algebra, when we associate 
bimodule $\mathbf{B}_i$ to the element $C_i'.$ Multiplication by 
$q$ corresponds to the grading shift up by $1.$ 

Furthermore, the Kazhdan-Lusztig basis in $H_3$ is given by 
\begin{eqnarray*} 
& & C_{e}' = 1, \hspace{0.1in} C_{s_1}' = C_1', \hspace{0.1in} 
 C_{s_2}' = C_2', \\
& & C_{s_1s_2}' = C_1' C_2', \hspace{0.1in} C_{s_2s_1}' = C_2' C_1', \\
& & C_{w}' = C_1' C_2' C_1' - C_1' = C_2'C_1'C_2' - C_2',
\end{eqnarray*} 
where $w=s_1s_2s_1=s_2s_1s_2$ is the longest element in $S_3$ and 
$e$ denotes the unit element of $S_3.$  

\vspace{0.1in} 

Arbitrary tensor products of bimodules $\mathbf{B}_1,\mathbf{B}_2$
can have only $6$ different indecomposable summands, up to isomorphism 
and grading shifts: 
$$ R, \hspace{0.1in} \mathbf{B}_1, \hspace{0.1in} \mathbf{B}_2,  
\hspace{0.1in} \mathbf{B}_1\otimes_R  \mathbf{B}_2, \hspace{0.1in}
\mathbf{B}_2\otimes_R  \mathbf{B}_1, \hspace{0.1in}
 \mathbf{B}_{1,2}.$$
Tensor products of these bimodules match the multiplication in 
the Kazhdan-Lusztig basis of $H_3.$ 

\vspace{0.1in} 

Soergel extends this patterns to all $m.$ For technical reasons he uses 
$\mathbb{C}$ as the ground field instead of $\Q.$ He shows the existence of 
graded indecomposable $R$-bimodules $\mathbf{B}_w,$ for $w\in S_m,$  
with the following properties: 
\begin{itemize} 
\item $\mathbf{B}_w$ is a finitely-generated projective left $R$-module
and a finitely-generated projective right $R$-module. 
\item $\mathbf{B}_{s_i}= \mathbf{B}_i,$ and $\mathbf{B}_{e}= R.$ 
\item This collection of bimodules is closed under tensor product:  
$$\mathbf{B}_{w}\otimes_R \mathbf{B}_y \cong \oplusop{z\in S_m} 
   \mathbf{B}_z^{n^z_{wy}},$$ 
where $n^z_{wy}\in \N[q,q^{-1}]$ are structure coefficients of 
the multiplication in the Kazhdan-Lusztig basis $\{ C_w'\}$ of $H_m,$
 $$C_w' C_y' =\sum_z n_{wy}^z C_z'.$$
\end{itemize} 

This construction, together with  the original work of Kazhdan and Lusztig, 
can be viewed as a categorification of the Hecke algebra. The 
Kazhdan-Lusztig basis lifts to a collection of bimodules, multiplication 
in the Hecke algebra lifts to the tensor products of bimodules, etc. 

More precisely, Soergel bimodules can be used to produce a categorification 
of the Hecke algebra action on its regular representation. The regular 
representation becomes the Grothendieck group of the graded version 
of the category of Harish-Chandra bimodules for $\mathfrak{sl}_m$ (with 
generalized trivial central character on both sides). To each $w\in S_m$ 
there is assigned an exact endofunctor in this category, which acts 
on the Grothendieck group in the same way as $C_w'$ acts by 
left multiplication on $H_m.$ The endofunctors can be reconstructed 
from bimodules $B_w,$ composition of endofunctors matching the 
tensor product of these bimodules. Some details can be found in [S1], others 
follow from Soergel's results.  

\hspace{0.1in} 

The group $G=SL(m,\C)$ acts transitively on the flag variety $G/B.$ 
The diagonal action of $G$ on $G/B\times G/B$ has finitely many 
orbits, which are in a natural bijection with elements of the symmetric 
group. The diagonal orbit corresponds to the trivial element, and 
the open orbit--to the maximal length permutation. Let $\mathcal{O}_w$ 
be the orbit associated with $w\in S_m.$ There exists a complex of 
sheaves $\mathrm{IC}(\overline{\mathcal{O}}_w)$ 
(the intersection cohomology sheaf), 
 supported on the closure of $\mathcal{O}_w.$ This complex of 
sheaves is $G$-equivariant and its stalk cohomology groups are constant 
along each orbit. According to Soergel [S2], bimodule $\mathbf{B}_w$  
is isomorphic to the $G$-equivariant cohomology of this complex,  
  $$ \mathbf{B}_w\cong \mathrm{H}_G(\mathrm{IC}
 (\overline{\mathcal{O}}_w)).$$ 
The bimodule structure comes from identifying $R$ with the $G$-equivariant 
cohomology of $G/B.$ 

\vspace{0.1in} 

\emph{Example:} when $m=2,$ the flag variety is $\mathbb{P}^1$ 
and $G$ acts on $\mathbb{P}^1\times \mathbb{P}^1$ with two 
orbits: the diagonal and its complement.  
The closure of each orbit is smooth, and the IC sheaf is the constant 
sheaf on the orbit's closure, shifted by the dimension of the orbit. 
Consequently, the equivariant cohomology groups of these IC sheaves 
are the equivariant cohomology groups of the diagonal and 
of $\mathbb{P}^1\times \mathbb{P}^1.$ They  are 
\begin{eqnarray*} 
& & \mathrm{H}_{G}(\mathbb{P}^1) \cong \mathrm{H}_B(\cdot) 
\cong \mathrm{H}_{SO(2)}(\cdot) \cong \C[y] \cong R, \\  & & 
 \mathrm{H}_{G}(\mathbb{P}^1\times \mathbb{P}^1) \cong 
\mathrm{H}_{B}(\mathbb{P}^1) \cong 
\mathrm{H}_{SO(2)}(\mathbb{P}^1)\cong \mathbf{B}_1, 
\end{eqnarray*} 
where dot denotes a point, and $y$ a generator of 
$\mathrm{H}^2(\mathbb{CP}^{\infty}, \Q).$ 
The group $SO(2)$ acts on $\mathbb{P}^1\cong S^2$ by rotations 
about the north pole-south pole axis. The action is free except at the 
poles, and the quotient by the action is naturally the interval $[-1,1].$ 
 The cohomology can be rewritten 
as $\mathrm{H}(S^2\times_{SO(2)}ESO(2)).$ 
The space $S^2\times_{SO(2)}ESO(2)$ maps onto 
the orbit space $[-1,1]$ of $S^2$ under the $SO(2)$ action. The 
fiber over each point other than $-1,1$ is contractible, while 
over $-1$ and $1$ the fiber is isomorphic to $\mathbb{CP}^{\infty}.$ 
Therefore, the equivariant cohomology is naturally isomorphic 
to the cohomology of the 1-point union of two copies of 
$ \mathbb{CP}^{\infty}.$ This cohomology can be identified 
with $R\otimes_{R_1} R$ as a graded $R$-bimodule, and even as 
a ring. 

\vspace{0.1in} 

For an arbitrary $m,$ the bimodule $\mathbf{B}_i$ is isomorphic to 
the $G$-equivariant cohomology groups of the ``thickened'' 
flag variety 
\begin{eqnarray*} 
& & Y=\{ (L_1, \dots, L_i, \dots , L_{m-1}, L_i') |  
\dim(L_j)=j, \dim(L_i')=i,  \\
& & 0 \subset L_1 \subset L_2 \subset \dots 
\subset L_{m-1}\subset \mathbb{C}^m, L_{i-1}\subset L_i'\subset L_{i+1}
  \}.
 \end{eqnarray*}

\vspace{0.2in} 

{\bf Bimodules $\mathbf{B}_w$ and link homology.} 

\noindent 
Our definition of link homology $\mathrm{HHH}(\sigma)$ used 
only tensor products of $\mathbf{B}_i$ rather than all $\mathbf{B}_w,$ 
which are the indecomposable summands of the tensor products. 
However, the Koszul resolution of indecomposable bimodule
$\mathbf{B}_{s_i s_{i+1}s_i}$ appear in  
the proof of the invariance of $\overline{H}(\sigma)$ under the 
third Reimeister move [KR2, Section 6]. 

Indecomposable bimodules $\mathbf{B}_w$ might prove useful in 
computations of link homology. Each term $F^j(\sigma)$ of the 
complex $F(\sigma)$ decomposes as a
direct sum of $\mathbf{B}_w,$ with various shifts and multiplicities, 
$$ F^j(\sigma) \cong \oplusop{w} \mathbf{B}_w^{n_w(\sigma,j)}, 
 \hspace{0.2in} n_w(\sigma,j) \in \N[q,q^{-1}].$$ 
Suppose we know these multiplicities $n_w(\sigma,j)$ and have 
the formula for the differential with respect to these direct sum 
 decompositions of $F^j(\sigma).$ For each $j$ the differential is described by 
an $n!\times n!$ matrix, with rows and columns enumerated by $w.$ 
The $(y,w)$-entry is itself a matrix with $n_y(\sigma,j+1)|_{q=1}$ rows 
and $n_w(\sigma,j)|_{q=1}$ columns describing the direct summand 
$$ \mathbf{B}_w^{n_w(\sigma,j)} \lra 
    \mathbf{B}_y^{n_y(\sigma,j+1)}$$
of the differential 
$$ \partial : F^j(\sigma) \lra F^{j+1}(\sigma).$$  
Each entry of the latter matrix is a homomorphism of bimodules 
$\mathbf{B}_w \lra \mathbf{B}_y$ of a particular degree. 

This horrendous complex can be simplified, by stripping off 
contractible summands of the form 
 $$ 0 \lra \mathbf{B}_w\{i\} \stackrel{1}{\lra} \mathbf{B}_w\{i\} 
   \lra 0.$$ 
Such a summand exists whenever there is an entry in the 
$j$-th matrix of matrices which is a nonzero complex multiple (recall we 
are working over $\C$) of the identity map 
of $\mathbf{B}_w\{i\}.$ 

\vspace{0.1in} 

Throwing out all contractible summands from $F(\sigma)$ results 
in a much smaller complex which we denote $F_{min}(\sigma).$ 
Up to isomorphism, $F_{min}(\sigma)$ does not depend on 
the order and choices of removed contractible summands. 
The reduction to $F_{min}(\sigma)$ is best done inductively 
on the length of $\sigma=\sigma_{i_1}^{\epsilon_1}\dots \sigma_{i_n}^{\epsilon_n}$
Once we found
 $F_{min}(\sigma_{i_1}^{\epsilon_1}\dots \sigma_{i_r}^{\epsilon_r}),$ 
for some $r<n,$ tensor it with $F(\sigma_{i_{r+1}}^{\epsilon_{r+1}})$ 
and reduce to minimal size. We start with $r=1$ and proceed until $r=n.$ 
The resulting minimal complex is isomorphic to $F_{min}(\sigma).$ 

\vspace{0.1in}

To determine $\mathrm{HHH}(\sigma),$ we take the Hochschild 
homology of each term $F^j_{min}(\sigma),$ arrange them 
into a complex and take its cohomology. 

\vspace{0.1in} 

\emph{Remark:} A similar algorithm to compute the sl(2) link homology was 
found by Dror Bar-Natan and implemented by him and 
Jeremy Green [BN1,2,3]. 
They represent a link as a composition of elementary tangles 
$L=t_1\dots t_n.$ In the language of [K], 
the invariant of a tangle $t_1\dots t_r$ is a complex 
of graded projective $H^s$-modules where $2s$ is the number of endpoints 
of $t_1\dots t_r.$ After splitting off all contractible summands from the complex
(thus reducing it to minimal size),
tensor it with the complex assigned to $t_{r+1},$ reduce the product to minimal 
size, and so on. We note that Bar-Natan and Green use a 
more refined and, at the same time, more geometric framework than that of 
rings $H^s.$ 

\vspace{0.1in} 

\emph{Example:} When $m=2$ and $\sigma=\sigma_1^n,$ the 
reduction leads to a simple computation of homology groups $\mathrm{HHH}(\sigma).$ 
The complex $F(\sigma)$ consists 
of $2^n$ terms which are tensor powers of $\mathbf{B}_1.$ 
The minimal complex has only $n+1$ indecomposable bimodules: 
$$ 0 \lra R\{2n\} \stackrel{d^0}{\lra} B_1\{2n-2\} \stackrel{d^1}{\lra} 
B_1 \{ 2n-4\} \stackrel{d^2}{\lra} \dots \stackrel{d^n}{\lra} B_1 \lra 0 $$ 
The differential is 
\begin{eqnarray*} 
   d^0(1) & = & 1\otimes y + y \otimes 1 , \\
   d^i (1\otimes 1) & = & 1\otimes y - y \otimes 1 \hspace{0.1in} \mathrm{odd} \hspace{0.05in} 
   i>0 , \\
    d^i (1\otimes 1) & = & 1\otimes y + y \otimes 1 \hspace{0.1in} \mathrm{even} 
   \hspace{0.05in}  i>0, 
\end{eqnarray*} 
where $R= \Q[y], y=x_1-x_2, R_1= \Q[y^2]$ and $B_1=\Q[y]\otimes_{\Q[y^2]}
 \Q[y].$ We take the Hochschild homology of each term, and the resulting 
complex splits into the direct sum of two complexes of graded $R$-modules 
(for $\mathrm{HH}_0$ and $\mathrm{HH}_1$). The complex of Hochschild 
homologies of zero degree is  
$$ 0 \lra R\{2n\} \stackrel{2y}{\lra} R\{2n-2\} \stackrel{0}{\lra}
R\{2n-4\}\stackrel{2y}{\lra} R\{2n-6\}    \stackrel{0}{\lra} \dots. $$ 
It has cohomology $\Q$ in bidegrees $(2n-2, 1), (2n-6, 3), \dots, (0,n),$ 
where the first grading corresponds to the degrees of variables $x_i$ and 
the second is cohomological. 

The Hochschild homology complex of degree one is  
$$ 0 \lra R\{2n+2\} \stackrel{1}{\lra} R\{2n+2\} \stackrel{0}{\lra}
R\{2n\}\stackrel{2y}{\lra} R\{2n-2\} \stackrel{0}{\lra} R\{2n-4\}
\stackrel{2y}{\lra} \dots. $$ 
   It has cohomology $\Q$ in bidegrees $(2n-2, 3), (2n-6, 5), \dots, (4, n).$ 

Hence, for the $(2,n)$-torus knot $L,$ homology $\mathrm{HHH}(L)$ 
has rank $n,$ as previously computed 
by Rasmussen (private communication). It was predicted in 
[GSV], [DGR] that the suitable HOMFLYPT homology groups of the $(2,n)$ torus
knot have rank $n.$ If $\overline{H}\cong \mathrm{HHH}$ 
is isomorphic to the theory conjectured to exist in [GSV] and [DGR], one 
would have an interesting link between perverse sheaves on flag varieties 
and the Gromov-Witten theory on the total space of the 
$\mathcal{O}(-1)\oplus  \mathcal{O}(-1)$ bundle over $\mathbb{P}^1.$ 

\vspace{0.2in} 

The approach to computing $\mathrm{HHH}(\sigma)$ 
via bimodules becomes significantly more challenging already for $m=3.$ 
With six indecomposable bimodules $\mathbf{B}_w$ 
the endomorphism algebra $\mathrm{End}_{R^e}(\oplusop{w\in S_3}
 \mathbf{B}_w)$
looks quite complicated. If we don't exclude trivial cases and use 
various symmetries, describing the multiplication   
$$ \mathrm{Hom}(\mathbf{B}_w, \mathbf{B}_y) \otimes 
   \mathrm{Hom}(\mathbf{B}_y, \mathbf{B}_z) \lra 
   \mathrm{Hom}(\mathbf{B}_w, \mathbf{B}_z)$$
directly requires dealing with $6^3=216$ possibilities for $w,y,z\in S_3.$ 
Also, knowing the endomorphism algebra is only the first (and not the most difficult) 
step in the algorithm. It's almost certain that describing $\overline{H}(\sigma)$ via 
the limit $N\to \infty$ of $SL(N)$ link homologies, as suggested by 
Gornik, Rassmussen and [DGR], 
together with Rasmussen's methods [Ra] for computing the latter, will get 
the job done much faster. 

\vspace{0.2in} 

{\bf Acknowledgements:} This paper is a write-up of the talk 
the author gave at the knot homology conference at UQAM in Montreal 
in September, 2005. I am grateful to Olivier Collin for 
organizing the conference and giving me an opportunity to present 
the above observations. Several ideas implicit or explicit in this paper 
(setting $a$ of [KR2] to $0,$ passing from the Koszul complex of a wide 
edge to its zeroth cohomology) had been independently suggested by 
Jacob Rasmussen. Theorem~\ref{thm} clarifies 
Maxim Vybornov's conjecture [V] that the link homology of [KR2] should 
be related to the category $\mathcal{O}$ and perverse sheaves 
on the flag variety. We were also influenced by Jozef Przytycki's 
 comparison of the Hochschild homology 
of truncated polynomial algebras and the link homology of $(2,n)$-torus 
knots [P]. The author would like to thank Charles Frohman, Sergei Gukov, 
Lev Rozansky and Ilya Shapiro for useful discussions. 
Partial support came from the NSF grant DMS-0407784.

\vspace{0.2in} 

{\bf References} 

\noindent
[BN1] D.~Bar-Natan, Khovanov's homology for tangles and cobordisms, 
\emph{Geom. Topol.}  9 (2005) 1443-1499, arxiv  math.GT/0410495. 

\noindent 
[BN2] D.~Bar-Natan,  I've computed Kh(T(9,5)) and I'm happy, 
talk at George Washington University, February 2005, available online. 

\noindent
[BN3] D.~Bar-Natan, Knot Atlas, http://katlas.math.toronto.edu/wiki. 

\noindent 
[BB] A.~Beilinson and J.~Bernstein, Localization de g-modules, \emph{C. R. 
Acad. Sci. Paris} (1) {\bf 292} (1981), 15--18.

\noindent
[BBD] A.~Beilinson, J.~Bernstein and P.~Deligne, Faisceaux 
pervers, \emph{Ast\'erisque} {\bf 100} (1982). 

\noindent 
[BK] J.~L.~Brylinski, M.~Kashiwara, Kazhdan-Lusztig conjecture 
and holonomic systems, \emph{Invent. Math.} {\bf 64} (1981), 387--410. 

\noindent
[D] P.~Deligne, Action du groupe des tresses sur une cat\'egorie, 
\emph{Inv. Math.} {\bf 128} (1997), 159--175. 

\noindent 
[DGR] N.~Dunfield, S.~Gukov and J.~Rasmussen, The superpolynomial for 
knot homologies, arxiv math.GT/0505662. 

\noindent
[GM] M.~Goresky and R.~MacPherson, Intersection homology II, 
\emph{Invent. Math.} 71, 77--129 (1983). 

\noindent
[GSV] S.~Gukov, A.~Schwarz and C.~Vafa, Khovanov-Rozansky homology 
and topological strings, arXiv hep-th/0412243. 

\noindent
[H] J.~Humphreys, Reflection groups and Coxeter groups, 
\emph{Cambridge studies in adv. math.} {\bf 29} 1990. 

\noindent 
[HOMFLY] P.~Freyd, D.~Yetter, J.~Hoste, W.~B.~R.~Lickorish, 
K.~Millett and A.~Ocneanu, A new polynomial invariant of knots 
and links, \emph{Bull. AMS. (N.S.)} {\bf 12} 2, 239--246, 1985. 

\noindent 
[KL1] D.~Kazhdan and G.~Lusztig, Representations of Coxeter groups 
and Hecke algebras, \emph{Invent. Math.} {\bf 53} (1980), 191--213. 

\noindent
[KL2] D.~Kazhdan and G.~Lusztig, Schubert varieties and Poincar\'e 
duality, \emph{Proc. Symp. Pure Math.} {\bf 36}, 1980.  

\noindent
[Ka] C.~Kassel, Homology and cohomology of associative 
algebras: a concise introduction to cyclic homology. Online notes. 

\noindent
[K] M.~Khovanov, A functor-valued invariant of tangles,
\emph{Algebr. Geom. Topol.} 2 (2002) 665-741,  math.QA/0103190. 

\noindent
[KR1] M.~Khovanov and L.~Rozansky, Matrix factorizations and
link homology, arxiv math.QA/0401268.

\noindent
[KR2] M.~Khovanov and L.~Rozansky, Matrix factorizations and
link homology II, arxiv math.QA/0505056. 

\noindent 
[KW] R.~Kiehl and R.~Weissauer, Weil conjectures, perverse sheaves 
and l'adic Fourier transform, Springer-Verlag, Berlin, 2001. 

\noindent 
[Ki] F. Kirwan, An introduction to intersection homology theory, 
\emph{Pitman Research Notes in Mathematics Series, vol. 187,} 
Longman Scientific and Technical, Harlow, 1988.

\noindent
[L] J.-L.~Loday, Cyclic homology, 2nd edition. Springer-Verlag, Berlin, 1998. 

\noindent
[M] R.~MacPherson, Global questions in the topology of 
singular spaces, \emph{Proc. Int. Cong. Math. Warszawa,} 1983. 

\noindent
[Mi] D.~Mili\v{c}i\'c, Localization and representation theory 
of reductive Lie groups,  http://www.math.utah.edu/ftp/u/ma/milicic/math/book.dvi

\noindent 
[P] J.~Przytycki, When the theories meet: Khovanov homology as 
Hochschild homology of links, arxiv math.GT/0509334. 

\noindent 
[PT] J.~Przytycki and P.~Traczyk, Conway Algebras and
Skein Equivalence of Links, \emph{ Proc. AMS}  {\bf 100}  744-748, 1987.

\noindent
[Ra] J.~Rasmussen, Khovanov-Rozansky homology of two-bridge knots 
and links, arxiv math.GT/0508510.

\noindent 
[R] R.~Rouquier, Categorification of the braid group, arxiv 
math.RT/0409593. 
 
\noindent
[Sa] M.~Saito, Mixed Hodge modules and applications, 
\emph{Proceedings of the International Congress of Mathematicians, 
vol. I, II (Kyoto, 1990),} 725--734, Math. Soc. Japan, Tokyo, 1991. 

\noindent 
[S1] W.~Soergel, The combinatorics of Harish-Chandra bimodules, 
\emph{J. reine angew. Math.} {\bf 429} (1992), 49--74. 

\noindent
[S2] W.~Soergel, Grading on representation categories, in: 
\emph{Proceedings of the ICM94 in Zürich,} Birkhäuser (1995), 800--806. 

\noindent
[Sp] T.A.~Springer, Queslques applications de la cohomologie 
d'intersection, \emph{S\'eminaire Bourbaki} 1981/82, 249--273. 

\noindent 
[T] T.~Tanisaki, Hodge modules, equivariant K-theory and Hecke 
algebras, \emph{Publ. RIMS, Kyoto Univ.} {\bf 23} (1987), 841--879. 

\noindent 
[V] M.~Vybornov, Private communication, June 2005. 

\noindent 
[W] H.~Wu, Braids and the Khovanov-Rozansky cohomology, 
arxiv math.GT 0508064. 

\vspace{0.2in} 

\noindent 
Mikhail Khovanov, Department of Mathematics, Columbia University, 
New York, NY 10027, khovanov@math.columbia.edu

\end{document}